\input amstex
\UseAMSsymbols
\input epsf.tex
\documentstyle{amsppt}


\loadbold \nologo \pageheight{8.5truein} \pagewidth{7.0truein}
\topmatter
\title $C_1$ in [2] is zero
\endtitle
\author Abbas Bahri
\endauthor
\endtopmatter

\medskip
\subhead{0.Introduction}\endsubhead
\medskip

This short Note falls within the context of the references [1], [2],
[3] and [4] cited below and establishes that the counter-example
provided in [4] to the second inequality of Corollary (19.10) in the
Clay Institute Monograph by J.Morgan and G.Tian [1] entitled "Ricci
Flow and the Poincare Conjecture" still stands after the correction
published in [2]. The constant $C_1$ in Lemma $0.4$ of [2] can be
taken to be zero.

The problem lies deeper.
\medskip
\subhead{1. Preliminaries}\endsubhead
\medskip

We assume in the sequel that the curve-shortening flow, starting
from a given curve, defines a piece of (immersed)surface $\Sigma$.
This happens for example  when $k(c(x_0,0))$, the norm of the
curve-shortening flow deformation vector $H(c(x,0))$ as in eg [1],
is non-zero at a given point $x_0$ of a smooth immersed curve
$c(x,0)$. Extending in section to the curve-shortening flow, we find
that an open set $U$ in $M$ is parameterized as $\{c_\mu(x,t)\}$,
$\mu$ an extra-parameter, with $\frac {\partial c_\mu(x,t)}{\partial
t}=H((c_\mu(x,t))=\nabla^{g(t)}_SS(c_\mu(x,t))$, $S$ is the unit
vector of $x \longrightarrow c_\mu(x,t)$, ($(t, \mu)$ frozen) for
the metric $g(t)$ evolving as in [1] through the Ricci flow.

$U$ is now mapped into $M\times [0, \epsilon)$ through the map
$c_\mu(x,t) \longrightarrow (c_\mu(x,t), t), t\in [0, \epsilon)$.

This is the framework of [2], with the metric $\hat {g}$ on $M
\times [0, \epsilon)$. The image of $M$ through this map will be
denoted $M_1$ in the sequel.

\medskip
\subhead {2. $C_1$ of [2] is zero}\endsubhead
\medskip

\subsubhead{a. The estimate $(\hat{\nabla}_{\frac
{\partial}{\partial
t}}H,H)_{(c(x,t),t)}=-Ric_{g(t)}(H,H)=O(k^2)$}\endsubsubhead
\medskip

 Consider the identity:

$$(\hat{\nabla}_{\hat H}H,H)=(\nabla^{g(t)}_H H,H)_{g(t)}-Ric_{g(t)}(H,H)$$

This identity is derived from the two ways to compute $\frac
{\partial k^2}{\partial t}$, the first one using the metric and
connection over $M \times [0, \epsilon)$, see [2], the second one by
differentiating directly as in (19.1) of [1] and using the
connection of the metric $g(t)$. The proof of (19.1) is repeated in
section 3, below, for the sake of completeness.

$\hat{\nabla}_{\hat H}H$ above is viewed as covariant
differentiation along the curve $(c(x,t),t)$ of $M_1$ (see section
1). Along this curve, $H=\nabla^{g(t)}_SS$, with $S=\frac {\frac
{\partial c(x,t)}{\partial x}}{|\frac {\partial c(x,t)}{\partial
x}|_{g(t)}}$. Since this quantity depends only on the value of $H$
along the curve $(c(x,t),t)$, $H$ can be extended in this subsection
to the $(c(x,t),s)$s as $H(c(x,t),s)=\nabla^{g(t)}_SS, S=\frac
{\frac {\partial c(x,t)}{\partial x}}{|\frac {\partial
c(x,t)}{\partial x}|_{g(t)}}$. This is understood as covariant
differentiation of $S$ along the curve parametrized in the variable
$x$ as $(c(x,t),s)$, with $t$ and $s$ frozen, hence as covariant
differentiation along $S$.

With $H=\nabla^{g(t)}_SS, S=\frac {\frac {\partial c(x,t)}{\partial
x}}{|\frac {\partial c(x,t)}{\partial x}|_{g(t)}}$, the expression
above $(\nabla^{g(t)}_H H,H)_{g(t)}$, at $t=t_0$, is in fact
$[(\nabla^{g(t_0)}_H H,H)_{g(t_0)}]_{t=t_0}$, see the proof below in
section 3.

Observe now that, since $H$ is horizontal:

$$\{[(\nabla^{g(t_0)}_H H,H)_{g(t_0)}](c(x,t))\}_{t=t_0}=(\hat{\nabla}_H H,H)_{(c(x,t_0),t_0)}$$

We use now the fact that $\hat H=\frac {\partial}{\partial t}+H$. We
can then split $\hat{\nabla}_{\hat H}$ into
$\hat{\nabla}_H+\hat{\nabla}_{\frac {\partial}{\partial t}}$ in the
derivation formula above. Using this splitting, we derive that

$$(\hat{\nabla}_{\frac
{\partial}{\partial t}}H,H)_{(c(x,t),t)}=-Ric_{g(t)}(H,H)=O(k^2)$$

This estimate can also be derived directly, without comparing the
two formulae, but with working on $M \times [0, \epsilon)$, with the
connection $\hat {\nabla}$ and a suitable coordinate frame, using
$\frac {\partial}{\partial t}, \frac {\partial c(x,t)}{\partial x},
H(c(x,t),s)=\nabla^{g(t)}_SS $, with $S(c(x,t),s)=\frac {\frac
{\partial c(x,t)}{\partial x}}{|\frac {\partial c(x,t)}{\partial
x}|_{g(t)}}$ and an additional suitable vector-field which adds the
parameter $\mu$ of section $1$.  We may assume that, with this
additional vector-field, we find a frame that is orthogonal (not
orthonormal) at $(c(x,t_0),t_0)$. The computation at $(c(x,t_0),
t_0)$ becomes straightforward: $H$ has a component on itself equal
to $1$; any derivative of this component is zero.  Then,
$(\hat{\nabla}_{\frac {\partial}{\partial
t}}H,H)_{(c(x,t_0),t_0)}=\hat {\Gamma}^{H}_{\frac
{\partial}{\partial t},H}(c(x,t_0),t_0)(H,H)_{g(t_0, c(x,t_0))}$.
The Lie bracket $[\frac {\partial}{\partial t}, H]$ is zero
\footnotemark \footnotetext {This follows from the action of the one
parameter group of $\frac {\partial}{\partial t}$; this action is
only on the second factor of the couple $(c(x,t),s)$, observe that
$H(c(x,t),s)=\nabla^{g(t)}_SS$, the value of $S$ being
$S(c(x,t),s)=\frac {\frac {\partial c(x,t)}{\partial x}}{|\frac
{\partial c(x,t)}{\partial x}|_{g(t)}}$, does not depend on $s$.
This is the framework, which we defined above, for our
computation.}; furthermore, the base point $c(x,t)$ for the
computation of the dot-products in $g(t)$ does not move under the
action of the one-parameter group of $\frac {\partial}{\partial t}$.
The only parameter that changes under the action of the one
parameter group of $\frac {\partial}{\partial t}$ is the second
parameter $s$ in $(c(x,t),s)$; over a time equal to $\tau$, this
involves a change of the metric from $g(t,c(x,t))$ into $g(t+\tau,
c(x,t))$. $H$ is unchanged. Thus, we find that as above:
$$(\hat{\nabla}_{\frac
{\partial}{\partial
t}}H,H)_{(c(x,t_0),t_0)}=-Ric_{g(t_0)}(H,H)=O(k^2)$$

Since $t_0$ is arbitrary, this new direct computation confirms the
former one.
 We have skipped some details in this construction and in this computation; they do not present any real difficulty.

Observe that the estimate above $O(k^2)$ does not depend on
$\Sigma$. It depends only on the ambient metric and on the value of
$k$ of course.

\bigskip
\subsubhead {b. $C_1$ in Lemma $0.4$ of [2] can be taken to be
zero}\endsubsubhead
\medskip

Coming back to the formula of the correction [2], we focus on the
addition of the two terms (observe the change of order in $S$ and
$H$ for the two last arguments in $\hat{Rm}$ with respect to [2] in
our notation. This is only a matter of notation: the content, which
is derived through the formula using the curvature operator
$(R(A,B)C,D)=Rm(A,B,C,D)$ in our notations and
$(R(A,B)C,D)=Rm(A,B,D,C)$ in the notations of [1] and [2] is
unchanged) $\hat {Rm}(\hat {H}, S, S, H)+(\hat{\nabla}_S
\hat{\nabla} _S \frac {\partial}{\partial t},H)$. These are the two
terms ( multiplied by $2$) contributing to the constant $C_1$ of
Lemma $0.4$ of [2]\footnotemark \footnotetext{ The other terms
displayed in [2] which could contribute to $C_1$ in Lemma $0.4$,
such as $2(\hat{\nabla}_{\hat{H}} [Ric_g(S,S)\frac
{\partial}{\partial t}], H)_g$, are in fact $O(k^2)$. For example,
for this latter term, we know that $H$ and $\frac
{\partial}{\partial t}$ are orthogonal and we know that
$\hat{\nabla}_{\frac {\partial}{\partial t}}(\frac
{\partial}{\partial t})$ is zero. $\hat{H}$ can be broken into
$\frac {\partial}{\partial t}+H$. The contribution of $\frac
{\partial}{\partial t}$ in the expression above is zero using our
observations and the contribution of $H$ is then $O(k^2)$ since
there is an additional dot product with $H$ in this expression.}.

  We need to be careful with
the definition of $S$: $S$, within the framework of the metric $\hat
{g}$ of [2], is equal to $S_2(c(x,t),s)= \frac {\frac {\partial
c(x,t)}{\partial x}}{|\frac {\partial c(x,t)}{\partial x}|_{g(s)}}$.
We will modify, without loss of generality for the computation, its
definition in the expression for $\hat {Rm}$ below.

Observe that the one parameter group of $\frac {\partial}{\partial
t}$, used over the time $\tau$, maps $(c(x,t),s)$ into $(c(x,t),
s+\tau)$. It follows that commutation of $\frac {\partial}{\partial
t}$ with $X=\frac {\partial c(x,t)}{\partial x}$ does occur, so that
$[\frac {\partial}{\partial t},X]=0$.

We now replace $S_2$, in $\hat{Rm}$ {\bf only}, by
$S(c(x,t),s)=\frac {\frac {\partial c(x,t)}{\partial x}}{|\frac
{\partial c(x,t)}{\partial x}|_{g(t)}}$, over $M\times [0,
\epsilon)$; since we are completing this computation at points of
$M_1$, where the value of $S_2$ is indeed $\frac {\frac {\partial
c(x,t)}{\partial x}}{|\frac {\partial c(x,t)}{\partial x}|_{g(t)}}$,
this modification is legitimate in $\hat{Rm}$: the value of  $\hat
{Rm}(\hat {H}, S_2, S_2, H)$ at a point of $M_1$ is the same than
the value of $\hat {Rm}(\hat {H}, S, S, H)$, $S_2$ and $S$ as above,
at the same point. This follows from the tensor properties of $\hat
{Rm}$.

As in [2], with new notations-$S_2$ is a notation which is not used
in [2]- $S$ in the other term $(\hat{\nabla}_S \hat{\nabla} _S \frac
{\partial}{\partial t},H)$ remains $S_2(c(x,t),s)= \frac {\frac
{\partial c(x,t)}{\partial x}}{|\frac {\partial c(x,t)}{\partial
x}|_{g(s)}}$. Observe that $S_2=\gamma S$, $\gamma=1, S_2.\gamma=0$
on $M_1$, so that $\hat {\nabla}_{S_2}\hat {\nabla}_{S_2}=\hat
{\nabla}_{S}\hat {\nabla}_{S}$ on $M_1$, see section 1, above for
the definition of $M_1$.

We need only to consider the terms which are not $O(k^2)$ and,
therefore, the above expression can be changed into $\hat {Rm}(\frac
{\partial}{\partial t}, S, S, H)+(\hat{\nabla}_{S_2} \hat{\nabla}
_{S_2} \frac {\partial}{\partial t},H)$.

 $H$ has been extended, without loss of generality, in both terms, over $M\times
[0,\epsilon)$ as $H(c(x,t),s)=\nabla^{g(t)}_SS$, the covariant
derivative along the unit vector $S$ for $g(t)$ to the curve $x
\longrightarrow (c(x,t),s)$, $(t,s)$ frozen.

 Writing the first term
with the use of the riemannian tensor on $\frac {\partial}{\partial
t}$ and $S$, this is $(\hat{\nabla}_{\frac {\partial}{\partial
t}}\hat{\nabla}_SS-\hat{\nabla}_S \hat{\nabla}_{\frac
{\partial}{\partial t}}S-\hat{\nabla} _{[\frac {\partial}{\partial
t},S]}S,H)$. This first term is thereby divided itself into three
further terms, which we now discuss one by one: We use the fact that
$[\frac {\partial}{\partial t},S]=\theta S$, with $\theta$ bounded.
This allows to see that the term $(\hat{\nabla} _{[\frac
{\partial}{\partial t},S]}S,H)$ is $O(k^2)$. The second term, after
the use of the commutation relation $[\frac {\partial}{\partial
t},S]=\theta S$ and the fact that $S$ and $H$ are orthogonal,
cancels with $(\hat{\nabla}_{S_2} \hat{\nabla} _{S_2} \frac
{\partial}{\partial t},H)$ (use our observation above about $\gamma,
S_2.\gamma$ on $M_1$) leaving $O(k^2)$.

Using the fact that $\hat{\nabla}_{\frac {\partial}{\partial
t}}\frac {\partial}{\partial t}=0$ and the fact that $\frac
{\partial}{\partial t}$ and $H$ are orthogonal, we find that the
first term is $(\hat{\nabla}_{\frac {\partial}{\partial t}}H,H)$,
which is, from our reasoning above, $O(k^2)$.

The conclusion follows. $C_1$ in [2] can be taken to be zero (we
assume that we have a piece of surface $\Sigma$; the estimate, as
pointed out above does not depend on $\Sigma$. When $k=0$ and there
is no immersed $\Sigma$ locally, the estimate of [2] contains $C_1
k$ and our argument does not work as is. It is however then clear
that $(\hat {\nabla}_{\frac {\partial }{\partial
t}}H,H)_{g(t)}=0=O(k^2)$ at such points. It is also clear that $C_1
k$ is zero at such points and can be forgotten. Using then density
and continuity, the assumption that $\Sigma$ immersed exists can be
removed). The additional terms in $k$ in the correction hide a
cancellation.

\bigskip
\subhead {3. Proof of (19.1) of [1]}\endsubhead
\medskip

For completion, we add here the computation that shows that (19.1)
of [1] holds: Let $\gamma(t)=g(t)-g(t_0)$ be the bilinear symmetric
2-tensor form defined by difference. We write on the image (piece
of) surface in M, which we assume to be immersed:

$$(H,H)_{g(t)}=(H,H)_{g(t_0)}+(H,H)_{\gamma(t)}$$

$H$ is here equal to $H=\nabla^{g(t)}_SS$, $S$ the unit vector
tangent to the curve $x \longrightarrow c(x,t)$, with respect to the
metric $g(t)$.

Differentiating with the use of the connection along the surface
induced by $g(t_0)$), we find, since $\frac {\partial}{\partial
t}=H$ on the surface:

$$[\frac {\partial k^2}{\partial t}]_{t=t_0}=[2(\nabla^{\Sigma}_HH,H)_{g(t_0)}+ \frac
{\partial ((H,H)_{\gamma(t)})}{\partial t}]_{t=t_0}$$

$\Sigma$ is our piece of surface.

Since $H$ is tangent to $\Sigma$, we may replace $\nabla^{\Sigma}$
with $\nabla^{g(t_0)}$ the connection for $g(t_0)$ on M. For the
derivative of the second term, we use local coordinates on $\Sigma$
and we find that this is $(H,H)_{\frac {\partial g}{\partial
t}_{t=t_0}}$, which yields the term $-2Ric(H,H)$. (19.1) follows.
Observe that, with $H_1=\nabla^{g(s)}_{S{_1}}S_1$, $S_1=\frac {\frac
{\partial c(x,s)}{\partial x}}{|\frac {\partial c(x,s)}{\partial
x}|_{g(s)}}$:

$$[(\nabla^{g(t_0)}_H\nabla_S^{g(t)}S,H)_{g(t_0)}]_{t=t_0}=[([\nabla^{g(t)}_HH_1]_{s=t},H)_{g(t)}]_{t=t_0}$$

\bigskip
\subhead {4. The Clay Institute computation [1], p442 for
$\frac{\partial k^2}{\partial t}$, slightly modified in order to
make it more transparent}\endsubhead
\medskip

The Clay Institute monograph [1], besides the division by $k$
pointed out in [3] and corrected in [2] (counter-example still
standing, see [4]), can be modified so that the computation of
$\nabla_HH$, p442, is more elementary and transparent. The final
result is unchanged: When computing $\nabla_HH$ as in p442 of [1],
it is preferable to take the connection $\nabla^{g(t)}$ at a fixed
value of $t$, $t=t_0$. This can be done, taking $H_2$ to be
$\nabla^{g(t_0)}_SS$, where $S$ is $\frac {\frac {\partial
c(x,t)}{\partial x}}{|\frac {\partial c(x,t)}{\partial x}|_{g(t)}}$,
The computation is carried as in [1], except that the Lie bracket
$[H_2,S]$ is not $\theta S$ anymore, since $H$ has been changed into
$H_2$. However, with $H_2$, $[H_2,S]= \theta S+O(t-t_0)$, where
$O(t-t_0)$ is small, as well as all its spatial derivatives as $t$
tends to $t_0$. There is no $t$-derivative in the computation of
[1]. The result for $\nabla^{g(t_0)}_{H_2}H_2$, p 442, reads ($H_2,
S$ as described above):

$$\nabla^{g(t_0)}_{H{_2}}H_2=\nabla^{g(t_0)}_S\nabla^{g(t_0)}_S
H_2+2(k^2+Ric(S,S))H+S.(k^2+Ric(S,S))S+Riem^{g(t_0)}(H,S)S+O(t-t_0)$$

In this computation now, $H_2$ can be replaced by $
H=\nabla^{g(t)}_SS$ in the left hand side. The result is unchanged
since there is no $t$-derivative in this formula, the corrective
terms are dropped into $O(t-t_0)$. $g(t_0)$ can be replaced by
$g(t)$ in the right hand side, where the covariant derivatives are
along $S$. Taking now the dot product with $H$, for the metric
$g(t_0)$ for the left hand side and for the metric $g(t)$ for the
right hand side, with $H=H_2$ there, we find
($(S,H_2)_{g(t)}=O(t-t_0)$):
$$(*) \quad (\nabla^{g(t_0)}_H H,H)_{g(t_0)}=(\nabla^{g(t)}_S\nabla^{g(t)}_S
H_2+2(k^2+Ric(S,S))H+Riem^{g(t)}(H,S)S+O(t-t_0),H_2)_{g(t)}$$

Understanding $\nabla^{g(t)}_S\nabla^{g(t)}_S H_2$ as covariant
differentiation along the curve $c(x,t)$, ie along $S$, we find that
$$(\nabla^{g(t)}_S\nabla^{g(t)}_S
H_2,H_2)_{g(t)}=(\nabla^{g(t)}_S\nabla^{g(t)}_S
H,H)_{g(t)}+O(t-t_0)$$

Then, after entering this and further replacing $H_2$ with $H$ in
$(*)$, taking then this modified identity at $t=t_0$, we find that
the computation of [1] for $\frac {\partial k^2}{\partial t}$ holds
at $t=t_0$. But $t_0$ is arbitrary.
\medskip
\subhead{5. Conclusion}\endsubhead
\medskip

In view of the present note and in view of the counter-example in
[4], the problem lies deeper.

\newpage

\widestnumber\no{99999}

\font\tt cmr12 at 24 truept 

\end{document}